\documentclass[10pt]{article}
\usepackage{amssymb}
\usepackage{latexsym}
\usepackage{amsthm}
\usepackage[dvips]{graphicx}
\usepackage{multicol}
\title{\textbf{A new approach to the logistic function with some applications} }
\author{ Grzegorz Rz\c{a}dkowski\thanks{Corresponding author: Email: g.k.rzadkowski@gmail.com}, Iwona G\l{}a\.zewska, Katarzyna Sawi\'nska 
}
\date{Warsaw University of Technology,\\Ludwika Narbutta 85, \\00-999 Warsaw, Poland} 
\begin{document}
\maketitle
\newtheorem{lemma}{Lemma}
\newtheorem{theorem}{Theorem}
\newtheorem{remark}{Remark}
\newtheorem{statement}{Statement}
\def \bangle{ \atopwithdelims \langle \rangle}  
\begin{abstract}
In the present paper we propose a new approach to investigate the logistic function, commonly used in mathematical models in economics and management. The approach is based on indicating in a given time series, having a logistic trend, some characteristic points corresponding to zeroes of successive derivatives of the logistic function. We give also examples of application of this method. 
\end{abstract}
\section{Introduction}The logistic equation is defined as
\begin{equation}\label{eq1}
	u'(t)=c_{1}u(u_{max}-u),\quad u(0)=u_{0}>0,
\end{equation}
where $t$ is time, $u=u(t)$ is an unknown function, $c_{1}, u_{max}$ are constants.  The constant $u_{max}$ is called the saturation level. The integral curve $u(t)$ fulfilling the condition $0<u(t)<u_{max}$ is known as the logistic function.\\
Many of economic phenomena, also related to the management follow the equation (\ref{eq1}) (see papers \cite{MI}, \cite{HCL}, \cite{MS}, \cite{QS}, \cite{WC}, \cite{Y}). \\
A phenomenon described by equation (\ref{eq1}) and a function $u(t)$ has an important property that its rate of growth $u'(t)$ is proportional to the level already achieved i.e. $u(t)$. On the other hand if $u(t)$ is sufficiently large then the factor $(u_{max}-u)$ is more and more significant and it influences inhibit further growth of the function $u(t)$. \\
Mathematically, equation  (\ref{eq1}) is the first order ordinary differential equation which is easily solved by a separation of variables method. \\
The main idea of the present paper, is to look, among the data of a given time series, for some characteristic points which correspond to zeroes of derivatives of the logistic function. One of these points is clearly the point corresponding to the inflection point (i.e. the zero of $u''$) of the logistic curve at which, as is well known, the logistic function takes the value $u_{max}/2$. For a sufficiently long time series the point corresponding to  the inflection point is easy to locate, even from the graph. If the data were collected for the time points spaced equally, then, instead of estimating the values of the first derivative, it is sufficient to calculate successive differences and seek the maximum for such a function.\\
 What we can do however, when the time series is not long enough, and we expect that the investigated phenomenon follows the logistic curve? When the phenomenon is on early stage of growth and the data is available only in a relatively short time interval? Statistical methods for estimating the parameters of the logistic function based, for example, on the method of the nonlinear least squares may be unreliable, since functions having significantly different values of the saturation level may produce slightly differing error values. A way to explanation of the situation seems in seeking, in the time series, points corresponding to zeroes of successive derivatives of the logistic function. For equally spaced time points this is equivalent to calculating successive differences. For example, as we will see in Sec. \ref{s3}, the zero of the third derivative $u'''$  (i.e., the extreme (maximum) of the second derivative $u''$) occurs at the point where the value of the logistic function is approximately $0.211\: u_{max}$.
\section{Logistic equation and logistic function}
\begin{figure}[h]
	\begin{center}
	 \includegraphics[height=4cm, width=7cm]{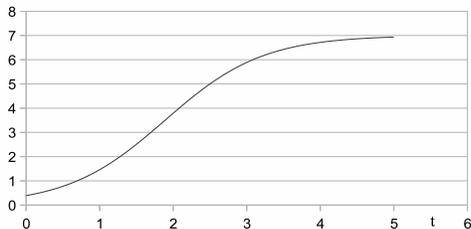}
	\end{center}
	\vspace{-12mm}
	\caption{Logistic function with parameters: $u_{max}=7,\; a=17$ and $c=1.5$}
	\label{p_0}
\end{figure}

We rewrite the logistic equation (\ref{eq1}) in the following, more convenient form, where the constant $c_{1}$, for computational reasons,  is written as $\displaystyle c_{1}=c/u_{max}$: 
\begin{equation}\label{eq2}
	u'(t)=\frac{c}{u_{max}}\:u(u_{max}-u),\quad u(0)=u_{0}.
\end{equation}
After solving (\ref{eq2}) we get the logistic function in the following form
\begin{equation}\label{eq3}
	u(t)=\frac{u_{max}}{1+ae^{-ct}},
\end{equation}
where constant $a$ appears in the integration process and is connected with the initial condition $\displaystyle u(0)=u_{0}=\frac{u_{max}}{1+a}$, therefore $\displaystyle a=\frac{u_{max}-u_{0}}{u_{0}}$. The Figure~\ref{p_0} shows the graph of an exemplary logistic function for $u_{max}=7,\; a=17$ and $c=1.5$.
\\In order to understand the further reasoning we have to introduce the so-called Eulerian numbers (see for instance Graham et al. \cite{GKP}). Let $\{a_1,a_2,\ldots ,a_n\}$ be a permutation of the set $\{1,2,\ldots ,n\}$. Then $\{a_{j},a_{j+1}\}$ is an ascent of the permutation if $a_j< a_{j+1}$. The Eulerian number denoted by $\displaystyle {n \bangle k} $ is defined as the number of permutations of the set $\{1,2,\ldots ,n\}$ having $k,\: (k=0,1,2,\ldots ,n-1)$ permutation ascents.  For example for $n=3$ the permutation $\{1,2,3\}$ has two ascents, namely $\{1,2\}$ and $\{2,3\}$, and  $\{3,2,1\}$ has no ascents. Each of the other four permutations of the set has exactly one ascent. Thus $\displaystyle { 3 \bangle 0 }=1 $,  $\displaystyle { 3 \bangle 1 }=4 $, and $\displaystyle { 3 \bangle 2 }=1 $. The first few Eulerian numbers are given in the Table~\ref{tab1}. It is well known that Eulerian numbers satisfy the following relations:
{\setlength\arraycolsep{2pt}
\begin{eqnarray*}
{ n \bangle k }&=&{ n \bangle n-k-1 }, \\
{ n+1 \bangle k }&=&(k+1){ n \bangle k }
+(n-k+1){ n \bangle k-1 },\\
{ n \bangle k } &=& \sum\limits_{j=0}^{k}(-1)^{j}{n+1 \choose j}(k-j+1)^{n}. 
\end{eqnarray*}} 
\begin{table}
\begin{center}
\begin{tabular}{|c|c c c c c c c c}
\hline
n & $\displaystyle {n \bangle 0} $ & $\displaystyle {n \bangle 1} $ & $\displaystyle {n \bangle 2} $ & 
$\displaystyle {n \bangle 3} $ & $\displaystyle {n \bangle 4} $ & $\displaystyle {n \bangle 5} $ & 
$\displaystyle {n \bangle 6} $ & $\displaystyle {n \bangle 7} $\\ \hline
0 & 1 & & & & & & &\\
1 & 1 & 0 & & & & & &\\
2 & 1 & 1 & 0 & & & & &\\
3 & 1 & 4 & 1 & 0  & & & &\\
4 & 1 & 11 & 11 & 1 & 0 & & &\\
5 & 1 & 26 & 66 & 26 & 1 & 0 & &\\
6 & 1 & 57 & 302 & 302 & 57 & 1 & 0 & \\ 
7 & 1 & 120 & 1191 & 2416 & 1191 & 120 & 1 & 0 \\
\hline
\end{tabular}
\end{center}
\caption{Eulerian numbers}
\label{tab1}
\end{table}
Equation (\ref{eq2}) is a particular case of Riccati's equation with constant coefficients
\begin{equation}\label{eq4}
	u'(t)=r(u-u_{1})(u-u_{2}).
\end{equation}
On the right hand side of (\ref{eq4}) is a quadratic function with the coefficient $r$ of $u^{2}$ and the roots $u_{1},\;u_{2}$. The constants $r\neq 0,\;u_{1},\;u_{2}$ can be generally the real or complex numbers. 
 
If $u(t)$ is a solution of (\ref{eq4}) then it is known a formula for the $n$th derivative $u^{(n)}(t)$ ($n=2,3,4,\ldots$)  of $u(t)$ expressing it in the function $u(t)$ itself:
{\setlength\arraycolsep{2pt}
\begin{eqnarray}
\hspace{-4mm}u^{(n)}(t) &=& r^{n}\left( { n \bangle 0 }
(u-u_{1})(u-u_{2})^{n}+{ n \bangle 1 }
(u-u_{1})^{2}(u-u_{2})^{n-1}\right. \nonumber \\
&& \left. +{ n \bangle 2 }
(u-u_{1})^{3}(u-u_{2})^{n-2}+\cdots 
+{ n \bangle n-1 }
(u-u_{1})^{n}(u-u_{2})\right)\nonumber \\
&=& r^{n}\sum\limits_{k=0}^{n-1}{ n \bangle k }
(u-u_{1})^{k+1}(u-u_{2})^{n-k}\label{eq5}
\end{eqnarray}}
where $n=2,3,\ldots $. \\
The above formula (\ref{eq5}) has been discussed during the Conference ICNAAM 2006
(September 2006) held in Greece and it appeared, with an inductive proof, in paper \cite{Rz1} (see also \cite{Rz2}). Independently the formula has been considered and proved, with a proof based on generating functions, by Franssens \cite{F}.\\
 The polynomial, of order $(n+1)$ of the variable $u$,  appearing on the right hand side of (\ref{eq5}) is known in the literature as the derivative polynomial. It can be proved (see \cite{Rz3}) that all $(n+1)$ roots of the polynomial are simple and lie in the interval $[u_{1}, u_{2}]$.  The derivative polynomials were recently intensively studied.
\section{Further properties of the logistic function and its derivatives}\label{s3}
\begin{figure}[h]
	\begin{center}
	 \includegraphics[height=4cm, width=7cm]{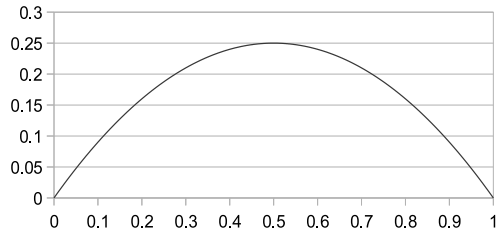}
	\end{center}
	\vspace{-10mm}
	\caption{Polynomial $P_{2}(u)$ for $u \in [0,\; 1]$}
	\label{p2}
\end{figure}

\begin{figure}[h]
	\begin{center}
	 \includegraphics[height=4cm, width=7cm]{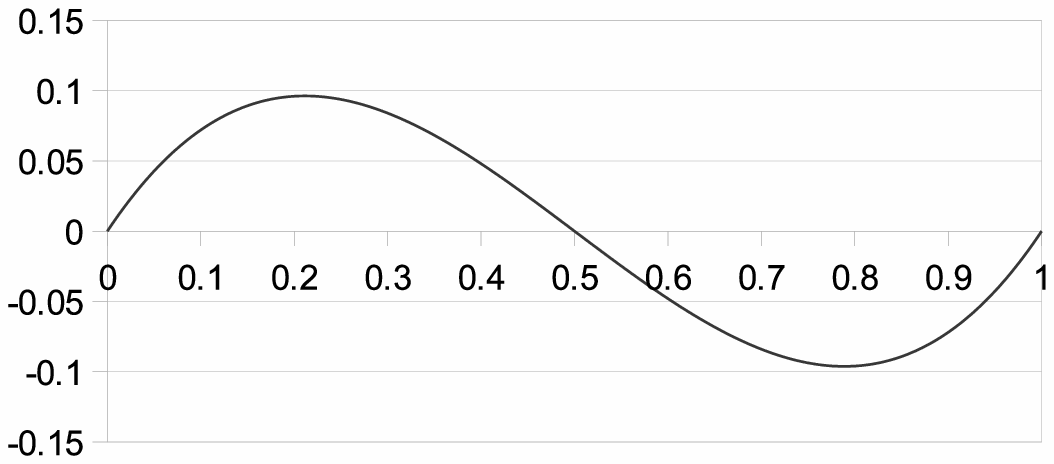}
	\end{center}
	\vspace{-12mm}
	\caption{Polynomial $P_{3}(u)$ for $u \in [0,\; 1]$}
	\label{p3}
\end{figure}

\begin{figure}[h]
	\begin{center}
	 \includegraphics[height=4cm, width=7cm]{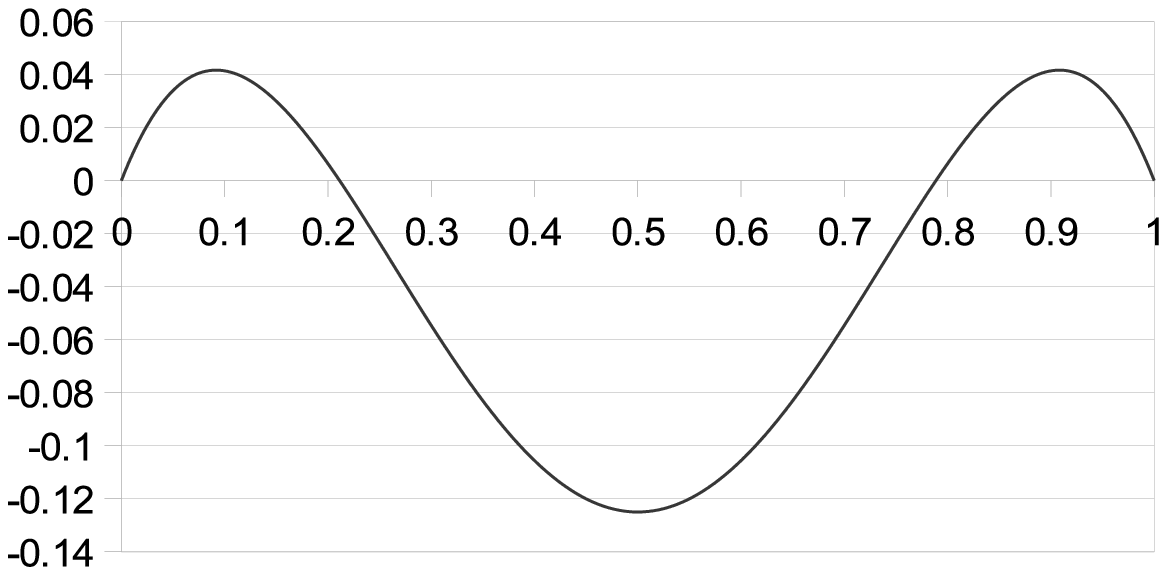}
	\end{center}
	\vspace{-10mm}
	\caption{Polynomial $P_{4}(u)$ for $u \in [0,\; 1]$}
	\label{p4}
\end{figure}

\begin{figure}[h]
	\begin{center}
	 \includegraphics[height=4cm, width=7cm]{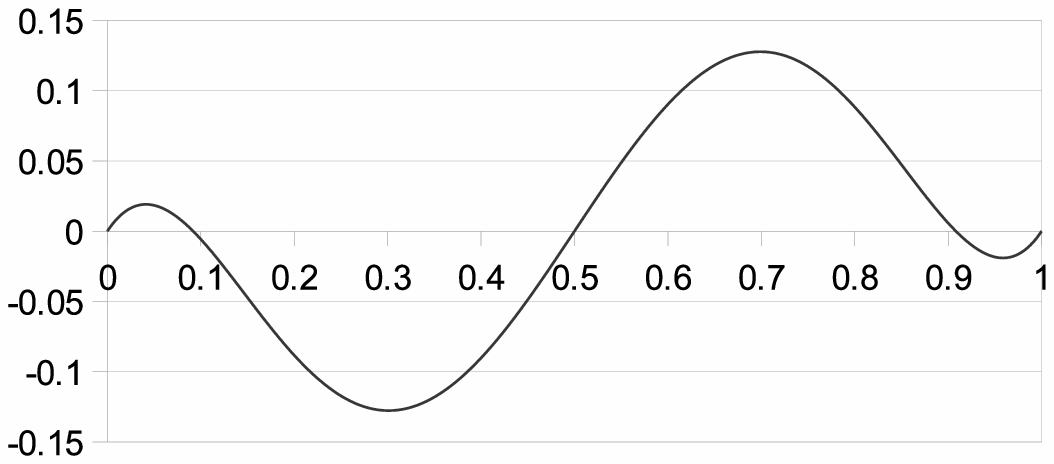}
	\end{center}
	\vspace{-12mm}
	\caption{Polynomial $P_{5}(u)$ for $u \in [0,\; 1]$}
	\label{p5}
\end{figure}

\begin{figure}[h]
	\begin{center}
	 \includegraphics[height=4cm, width=7cm]{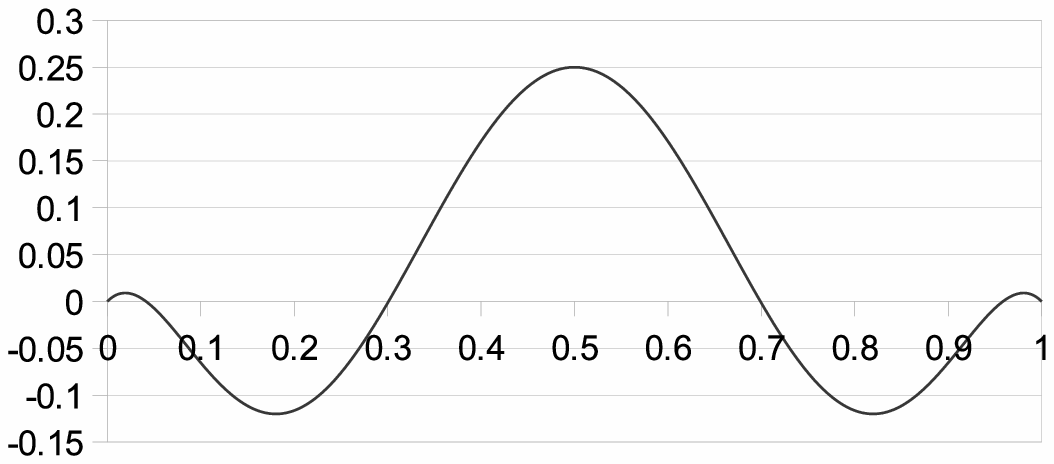}
	\end{center}
	\vspace{-12mm}
	\caption{Polynomial $P_{6}(u)$ for $u \in [0,\; 1]$}
	\label{p6}
\end{figure}

\begin{figure}[h]
\begin{center}
\includegraphics[height=4cm, width=7cm]{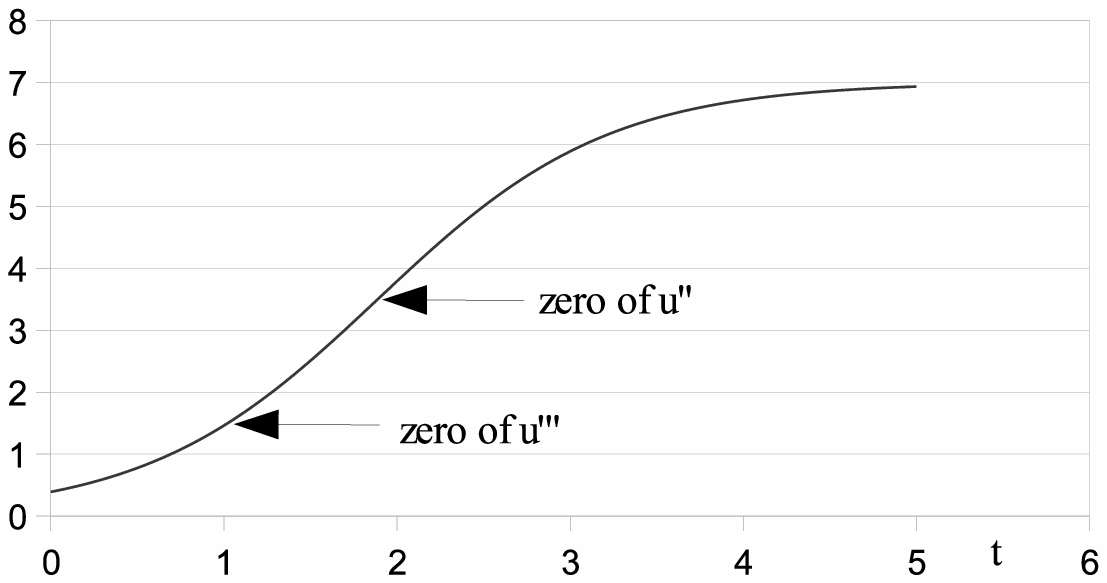}
\end{center}
\vspace{-12mm}
\caption{Two characteristic points of the logistic curve}
\label{p_00}
\end{figure}

Formula (\ref{eq5}) applied to the particular case of the logistic equation (\ref{eq2}) is as follows:
\begin{equation}\label{eq6}
	u^{(n)}(t)= \left(-\frac{c}{u_{max}} \right)^{n}\;\cdot\sum\limits_{k=0}^{n-1}{ n \bangle k }
u^{k+1}(u-u_{max})^{n-k}.
\end{equation}
The polynomial of the variable $u$ and of order $(n+1)$ on the right hand side of (\ref{eq6}) is uniform in the sense of the following.
\begin{remark}
If a number $u_0$ is a root of the polynomial on the right hand side of (\ref{eq6}) i.e.
\begin{equation}\label{eq7}
	\sum\limits_{k=0}^{n-1}{ n \bangle k }
u_0^{k+1}(u_0-u_{max})^{n-k}=0,
\end{equation}
then dividing both sides of (\ref{eq7}) by $u_{max}^{n+1}$ we get
	\[\sum\limits_{k=0}^{n-1}{ n \bangle k }
\left(\frac{u_0}{u_{max}}\right)^{k+1}\left(\frac{u_0}{u_{max}}-1\right)^{n-k}=0.
\]
Thus $u_0$ is a root of the derivative polynomial on the right hand side of (\ref{eq6}) if  $u_0/u_{max}$ is the root of the polynomial
\begin{equation}\label{eq8}
	P_{n+1}(u):=(-1)^{n}\sum\limits_{k=0}^{n-1}{ n \bangle k }
u^{k+1}(u-1)^{n-k}.
\end{equation}
\end{remark}
Let us write down, using formula (\ref{eq6}) and the notation of (\ref{eq8}), the first few derivatives of the logistic function, which fulfills equation (\ref{eq2}). By Remark 1 we can assume, without loss of the generality, that $u_{max}=1$ and $c=1$. \\
We obtain successively:
{\setlength\arraycolsep{2pt}
\begin{eqnarray*}
u'(t)&=&u(1-u)=-u(u-1)=P_{2}(u),\\[2mm]
u''(t)&=&u(u-1)^{2}+u^{2}(u-1)=P_{3}(u),\\[2mm]
u'''(t)&=&-u(u-1)^{3}-4u^{2}(u-1)^{2}-u^{3}(u-1)=P_{4}(u),\\[2mm]
u^{(4)}(t)&=&u(u-1)^{4}+11u^{2}(u-1)^{3}+11u^{3}(u-1)^{2}+u^{4}(u-1)=P_{5}(u),\\[2mm]
u^{(5)}(t)&=&-u(u-1)^{5}-26u^{2}(u-1)^{4}-66u^{3}(u-1)^{3}-26u^{4}(u-1)^{2}\\
&&-u^{5}(u-1)=P_{6}(u).
\end{eqnarray*}}
All roots of the polynomials $P_{k}(u) \;\textrm{for}\; k=3,4,5,6\;$ can be calculated explicitely, so the polynomials can be factored and we get  
{\setlength\arraycolsep{2pt}
\begin{eqnarray*}
P_{3}(u)&=& 2u(u-1)\left(u-\frac{1}{2}\right),\\[2mm]
P_{4}(u)&=&-6u(u-1)\left(u-\frac{1}{2}-\frac{\sqrt{3}}{6}\right)\left(u-\frac{1}{2}+\frac{\sqrt{3}}{6}\right),\\[2mm]
P_{5}(u)&=&24u(u-1)\left(u-\frac{1}{2}\right)\left(u-\frac{1}{2}-\frac{\sqrt{6}}{6}\right)\left(u-\frac{1}{2}+\frac{\sqrt{6}}{6}\right),\\[2mm]
P_{6}(u)&=&-120u(u-1)\!\!\left(u-\frac{1}{2}-\frac{\sqrt{30(15-\sqrt{105})}}{60}\right)\!\!\left(u-\frac{1}{2}-\frac{\sqrt{30(15+\sqrt{105})}}{60}\right) \\
&&\left(u-\frac{1}{2}+\frac{\sqrt{30(15-\sqrt{105})}}{60}\right)\!\!\left(u-\frac{1}{2}+\frac{\sqrt{30(15+\sqrt{105})}}{60}\right).
\end{eqnarray*}}
Therefore the least \emph{positive} root of the polynomial
{\setlength\arraycolsep{2pt}
\begin{eqnarray*}
P_{4}(u)&\textrm{is}&\;\; \frac{1}{2}-\frac{\sqrt{3}}{6}\approx 0.211,\\[2mm]
P_{5}(u)&\textrm{is}&\;\; \frac{1}{2}-\frac{\sqrt{6}}{6}\approx 0.0917,\\[2mm]
P_{6}(u)&\textrm{is}&\; \; \frac{1}{2}-\frac{\sqrt{30(15+\sqrt{105})}}{60} \approx 0.0413.
\end{eqnarray*}}
Thus by using Remark 1 we see for example that if at some point of time $t_0$ (the least possible) $u'''(t_0)=0$ ($u''(t_0)$ is a local maximum) then the value of the logistic function at this point is $u(t_0)=0.211\:u_{max} $. In Figure~\ref{p_00} we see two characteristic points of the exemplary logistic curve (with the same parameters as on Figure~\ref{p_0}): the inflection point (the zero of the second derivative $u''(t)$) and the zero of the third derivative $u'''(t)$.  Similar conclusions can be drawn for the least zeroes of the $u^{(4)}(t)$  (the polynomial $P_{5}(u)$ is used in this case) or $u^{(5)}(t)$ ($P_{6}(u)$) with the constants given above.\\
In Figures~\ref{p2}--\ref{p6} below we can see graphs of the derivative polynomials  $P_{k}(u)$ for $k=2,3,4,5,6$ respectively for $u \in [0,\; 1]$. The polynomials are symmetric (even or odd) with respect to the point $u=1/2$.
\section{Some applications of the method}
\begin{figure}
	\begin{center}
	 \includegraphics[height=4cm, width=7cm]{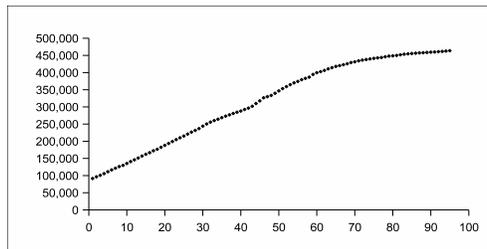}
	\end{center}
	\vspace{-12mm}
	\caption{The total number of loyalty cards issued (TNLC)}
	\label{p7}
\end{figure}

\begin{figure}
	\begin{center}
	 \includegraphics[height=4cm, width=7cm]{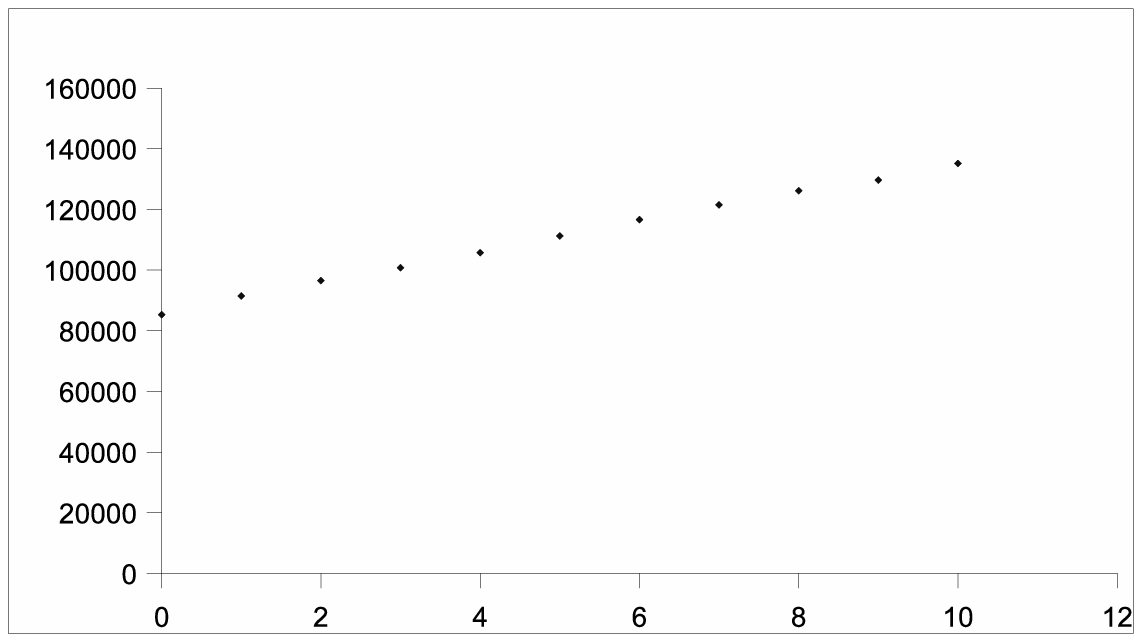}
	\end{center}
	\vspace{-12mm}
	\caption{Initial data }
	\label{p8}
\end{figure}

\textbf{\large Loyalty cards in a chain of stores}\\
The data in the Table~\ref{tab2}  represent the number of loyalty cards (NLC) issued in a large chain of stores in Poland. The observations relate to the period December 2011 - November 2013 (e.g., 47/12 means forty-seventh week of the year 2012). 
In the initial period of time, covering the first two or three months, the phenomenon was a fairly rapid process, due to Christmas, New Year and a big promotional campaign. Then the situation stabilized and in the next weeks the total number of  loyalty cards (TNLC) proceeded according to a logistic curve. Figure~\ref{p7} shows the total number of issued cards starting from the tenth week of all observations.
\\We will try to use, to the above data, our theory described in the earlier sections.  Let us take into account the first ten observations from the Figure~\ref{p7} (see the second and third columns of the Table 3 and Figure~\ref{p8}). We calculate the successive differences for the data. By the  second central difference (SCD) at time $t$ of an equally spaced time series $\{y_{t}\}$, ($t=0,1,2,3,\ldots,n$) at a time $t$ we mean the number given by the formula $(y_{t+1}-2y_{t}+y_{t-1})/2$.  The calculations are shown in Table~\ref{tab3}. Instead of using CSD we could use the  second left differences (SLD) given at time $t$ by the formula $(y_{t}-2y_{t-1}+y_{t-2})/2$.
\\We see that the first local maximal value of the second central difference $358$, is taken for $t=3$ where the total number of issued loyalty cards is $100,776$. Therefore, using comments from Section 3, we can estimate the saturation level as $u_{max}=100,776/0.211=477,611$. Since the value of SCD at point $t$ is equal to the value of SLD at $t+1$, then using in the decision the last one, leads to the estimation  $u_{max}=105,729/0.211=501,085$.
\\Otherwise we could find a polynomial, which best fits the data in the sense of the LSM and then investigate its second derivative for a maximum. Such polynomial e.g., of order four,  is as follows
	$$f(x) = -1.6807x^{4} + 30.515x^{3} - 206.23x^{2} + 5610.1x + 85584$$
and its second derivative $f''(x)=-20.16x^{2}+183.09x-412.46$ has a maximum at the point $x=4.5$. The value of $f(x)$ at this point is $f(4.5)=108,745$.  Thus we can estimate the saturation level of this phenomenon as $u_{max}=108,745/0.211=514,647$. \\
Starting a project to introduce a loyalty card e.g., in a chain of stores, we could use the logistic curve to determine the increase in the number of loyal customers.
At the same time, we could develop a promotion plan for the recruited persons, so that the budget of a given period of time (e.g., for a year) had a chance to be in a real way achieved.
During the project, the model created by the logistic curve allows us to monitor the effectiveness of the project, to draw conclusions and make appropriate decisions if there are derogations from our previous assumptions.\vspace{3mm}

\begin{figure}
\begin{center}
\includegraphics[height=4cm, width=7cm]{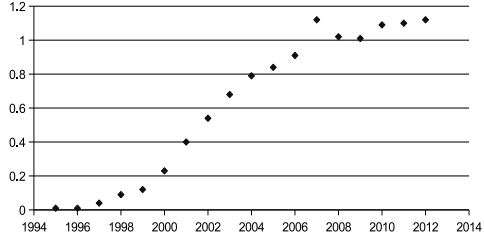}
\end{center}
\vspace{-9mm}
\caption{Rate of mobile telephone subscriptions per 1 inhabitant in Slovak Republic}
\label{p9}
\end{figure}

\begin{figure}
\begin{center}
\includegraphics[height=4cm, width=7cm]{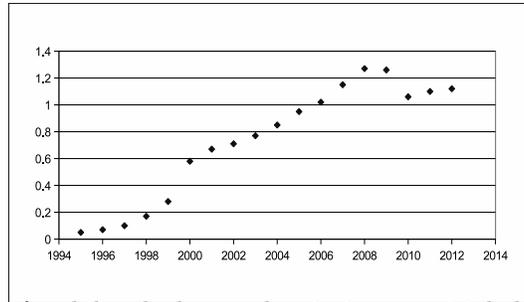}
\end{center}
\vspace{-9mm}
\caption{Rate of mobile telephone subscriptions per 1 inhabitant in Germany}
\label{p10}
\end{figure}

\begin{figure}
\begin{center}
\includegraphics[height=4cm, width=7cm]{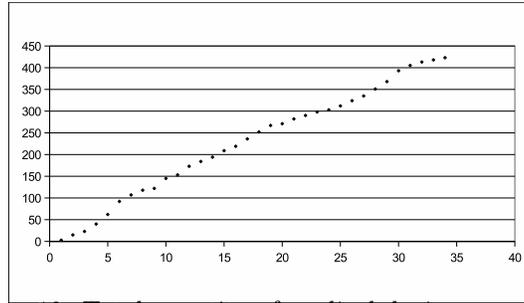}
\end{center}
\vspace{-9mm}
\caption{Total quantity of medical devices purchased}
\label{p11}
\end{figure}

\textbf{\large Diffusion of mobile telephones in two countries}\\
Recently  many papers have been published, devoted to mathematical modeling of the percentage diffusion over the population of mobile telephony for different countries (see \cite{HCL},\cite{MS},\cite{WC}, \cite{Y}). Consider two European countries, one with a stable and well-developed economy - Germany and the second, relatively recently accepted into the European Union -Slovak Republic, which transformed from a centrally planned economy to a market-driven economy. Table~\ref{tab5} shows the rate of mobile telephone subscriptions per 1 inhabitant in the two countries (see also Figures~\ref{p9} and \ref{p10}). The data, corresponding to the period
from 2000 to 2012, were collected from the International Telecommunication Union (ITU, http://www.itu.int) and these  corresponding to the period 1995-1999 were extracted from the paper C. Michalakelis, T. Sphicopoulos \cite{MS}. \\
On the graph for Germany (Figure~\ref{p10}) are seen two perturbations and changes in the trend. They were probably caused by early 2000s recession, which mainly occurred in developed countries and financial crisis of 2007–-2008, which led to the 2008–-2012 global recession. However in years 1995--2000 the shape well fits the logistic curve and we see from Table~\ref{tab5} that SCD clearly takes its maximal value at the point $0.28$. Thus the estimated saturation level is $u_{max}=0.28/0.211=1.327$. It seems interesting to note that, despite of the perturbations described above, a similar level has been reached in 2013, and was equal $1,301$\footnote{"Research and Markets Adds Report: Germany — Telecoms, IP Networks and Digital Media". TMC News, 12 June 2013. Retrieved 5 November 2013.}.\\
For fast-growing economy of Slovak  Republic the crises have not had much impact on the level of diffusion. Maximum od SCD, for the initial observations, is not as explicit as previously and should be fixed somewhere in the interval $[0.12, 0.23]$. Therefore the estimated level of saturation is from $0.12/0.211=0.57$ to $0.23/0.211=0.109$. Let us note that using SLD indicator would give us, in this case, the value of saturation level in the interval from $0.23/0.211=1.09$ to $0.4/0.211=1.89$.
\vspace{3mm} \\
\textbf{\large Purchases of certain medical devices}\\
The data in Table~\ref{tab6} relate to some specific medical devices used in the diagnosis and treatment of patients. These products are used in public and private medical institutions. The main users are public institutions, which buy these products using public funds in accordance with the law "Public Procurement Law". The average lifetime of the product is about three years. After this period, the device is subjected to a major renovation restoring its full functionality or is exchanged for a new one. The specifity of purchases from the budget indicates that purchases of these products are usually made in the fourth quarter of the calendar year. Increasing demand for these products contribute to periodic health programs implemented under the national program for the eradication of cancer.\\
Figure~\ref{p11} shows the total quantity of purchases (on the horizontal axes the first month denotes 06/09 i.e., June 2009).  We see that only the initial data fit a logistic curve and then the phenomenon is no longer of such a nature. However the maximum od SLD is reached in November 2009 where the total purchases are 92 devices. The estimated saturation level is $u_{max}=92/0.211=436$ and it seems to be a quite good forecast of its real value.
\section{Conclusions and future works} 
In this paper we have presented a method to study time series having the logistic trend. The approach is based on indicating in the data some characteristic points corresponding to zeroes of successive derivatives of the logistic function. We gave a general formula for the $n$th derivative of the logistic function, expressing it in terms of the derivative polynomials and Eulerian numbers. Then we have computed first few derivatives and calculated zeroes of the corresponding derivative polynomials. It seems that from a practical point of view, especially with relatively small number of initial data, particularly significant is the leading zero of the third derivative, where the logistic funtion takes value $0.211u_{max}$. 

We have shown the usefulness of our method with examples related to economics and management. We have demonstrated that when a phenomenon is at its early stages of development, then the saturation level may be effectively predicted. 
We believe that this approach should be used together with existing methods, for example, the nonlinear least squares method.\\
In our next research works we will apply a similar idea for other mathematical models used in economics and management as Gompertz and Bass curves.
\vspace{3mm} \\
\textbf{Conflict of Interests} The authors declare that they do not have any conflict of interest in their submitted
manuscript.

\begin{table}
\begin{center}
\begin{tabular}{||c|r|| c| r|| c |r ||c|  r ||}
\hline
week & NLC & week & NLC & 
week &  NLC & week & 
 NLC\\ \hline
48/11 & 7236 & 23/12& 4307&50/12 &7741 &25/13 &2077 \\
49/11 & 11904 & 24/12 &5776 & 51/12& 8950& 26/13&1889	\\
 50/11& 12887 & 25/12 & 5561 & 52/12&3447 &27/13 & 1686\\
51/11 & 10665 & 26/12 &5521 & 01/13  & 3510	&28/13 &1651 \\
 52/11& 5616 & 27/12 &5525 & 02/13 &6334  & 29/13&1402 \\
01/12 &  7133	& 28/12 &5625	& 03/13 &6793  & 30/13 &1247 \\
02/12 &8428  & 29/12 & 5393& 04/13 & 6846 & 31/13 & 2026\\ 
03/12 & 7263 & 30/12 & 5132 & 05/13 &5764  & 32/13 & 1847  \\
04/12 &7135&31/12&5768&06/13&5803	&33/13&899\\
05/12 &7038&32/12&5826&07/13&5121&34/13&1132\\
06/12 &6173	&33/12&4683	&08/13&4223&35/13&1920\\
07/12 &5061&34/12&5337&09/13&4955&36/13&1551\\
08/12 &4237&35/12&7216&10/13&3939&37/13&1172\\
09/12 &4953&36/12&6396&11/13&3566	&38/13&935\\
10/12 &5536&37/12&5325&12/13&7844&39/13&903\\
11/12 &5387	&38/12&4421	&13/13&5085&40/13&826\\
12/12 &4868&39/12&4111&14/13&3158&41/13&619	\\
13/12 &4673&40/12&4343&15/13&3550&42/13&840\\
14/12 &3496&41/12&4462&16/13&4468	&43/13&701\\
15/12 &5474&42/12&3780&17/13&3498&44/13&601\\
16/12 & 5576 &43/12&4048 &18/13&3726&45/13&882\\
17/12 &5245&44/12&3708&19/13&2339&46/13&775	\\
18/12 &5196&45/12&3474&20/13&2628&47/13&849\\
19/12 &5563&46/12&4462&21/13&2708	&48/13&1238\\
20/12 &5252&47/12&3957&22/13&3482&&\\
21/12 &4616	&48/12&5405	&23/13&2142&&\\
22/12 &5690&49/12&7913&24/13&2710&&\\
\hline
\end{tabular}
\end{center}
\caption{The number of issued loyalty cards (NLC)}
\label{tab2}
\end{table}
\begin{table}
\begin{center}
\begin{tabular}{|c|c |r |r| r |}
\hline
Week & No (t) &  TNLC &  SCD
\\ \hline
05/12 & 0 & 85305 &  \\
06/12& 1 & 91478 &  -556\\
07/12 & 2 & 96539 &  -412 \\
08/12 & 3 & 100776 &  358  \\
09/12 & 4 & 105729	&   291.5\\
10/12 & 5 & 111265 &  -74.5	 \\
11/12 & 6 & 116652	&  -259.5\\ 
12/12 & 7 & 121520	& -9.5  \\
13/12& 8& 126193	&  -588.5 \\
14/12& 9&129689	&  \\
\hline
\end{tabular}
\end{center}
\caption{Central second differences for the initial observations}
\label{tab3}
\end{table}
\begin{table}
\begin{center}
\begin{tabular}{||c||c |c ||c| c ||}
\hline
Year & Germany & SCD  & Slovak Republic & SCD
\\ \hline
1995 & 0.05&  &0.01 & \\
1996& 0.07 & 0.005 &0.01 &0.015 \\
1997 &  0.1& 0.02 & 0.04 &  0.01\\
1998 & 0.17 &  0.02& 0.09	& -0.01  \\
1999 & 0.28 & 0.095	& 0.12&0.04  \\
2000 & 0.58 & -0.105 & 0.23 &	0.03 \\
2001 & 0.67 & -0.025	&0.4 & -0.015 \\ 
2002 & 0.71 & 0.01	& 0.54& 0  \\
2003& 0.77& 0.01	&0.68 & -0.015 \\
2004&0.85 &0.01	&0.79 &  -0.03 \\
2005& 0.95&-0.015 & 0.84&0.01  \\
2006& 1.02& 0.03&0.91 & 0.07 \\
2007&1.15 & -0.005& 1.12&-0.155  \\
2008& 1.27& -0.065& 1.02& 0.045 \\
2009&1.26 &-0.095 &1.01 & 0.045 \\
2010& 1.06& 0.12& 1.09& -0.035 \\
2011& 1.1& -0.01&1.1 & 0.005 \\
2012&1.12 & &1.12 &  \\
\hline
\end{tabular}
\end{center}
\caption{Diffusion of mobile telephony for Germany and Slovak Republic}
\label{tab5}
\end{table}
\begin{table}
\begin{center}
\begin{tabular}{||c|r|| c| r|| c |r ||c|  r ||}
\hline
month & QMD & month & QMD & 
month &  QMD & month & 
 QMD\\ \hline
06/09 & 3 & 03/10& 23&12/10 & 15&09/11 &16 \\
07/09 & 12 & 04/10 &8 & 01/11& 4& 10/11&	17\\
 08/09&  8& 05/10 & 20 & 02/11&11 &11/11 & 25\\
09/09 & 17 & 06/10 & 11& 03/11  & 8	&12/11&12 \\
 10/09& 22 & 07/10 & 10& 04/11 & 8 & 01/12&8 \\
11/09 &  30	& 08/10 &15	& 05/11 &5  & 02/12 &5 \\
12/09 & 15 & 09/10 &10 & 06/11 &9  & 03/12 & 5\\ 
01/10 & 11 & 10/10 &17  & 07/11 & 12 &  &   \\
02/10 & 4 & 11/10 & 16 & 08/11 & 11 &  &   \\
\hline
\end{tabular}
\end{center}
\caption{The quantity of  medical devices purchased (QMD)}
\label{tab6}
\end{table}
\end{document}